\documentclass[12pt]{article}
\usepackage{amsfonts}
\usepackage{amssymb}
\usepackage{amsmath}
\usepackage{eufrak}
\usepackage[dvips]{graphicx}
\usepackage{latexsym}
\newcommand{\C}{{\EuFrak C}}

\newcommand{\bi}{\begin{itemize}}
\newcommand{\ei}{\end{itemize}}

\newtheorem{theorem}{Theorem}[section]
\newtheorem{lemma}[theorem]{Lemma}
\newtheorem{fact}[theorem]{Fact}

\newtheorem{corollary}[theorem]{Corollary}
\newtheorem{proposition}[theorem]{Proposition}
\newtheorem{definition}[theorem]{Definition}
\newtheorem{remark}[theorem]{Remark}

\newtheorem{question}[theorem]{Question}

\title{On $\omega$-categorical, generically stable groups}
\author{Jan Dobrowolski and Krzysztof Krupi\'nski\footnote{Research supported by the Polish Government grant N N201 545938}}
\date{}

\begin{document}
\maketitle
\begin{abstract}
We prove that each $\omega$-categorical, generically stable group  is solvable-by-finite. 
\end{abstract}
\footnotetext{2010 Mathematics Subject Classification: 03C45, 20A15}
\footnotetext{Key words and phrases: $\omega$-categorical group, generically stable type}

\section{Introduction}

A general motivation for us is to understand the structure of $\omega$-categorical groups satisfying various natural model-theoretic assumptions. There is, of course,  a long history of results of this kind, for example in a superstable, supersimple or NSOP context \cite{1,EW,6,9}  (see \cite[Introduction]{5} for a very quick overview of such results).

Recall that each countable, $\omega$-categorical group has a finite series of characteristic (i.e. invariant under the automorphism group) subgroups in which all successive quotients are characteristically simple groups (i.e. they do not have non-trivial, proper characteristic subgroups). On the other hand, Wilson (see \cite{10}) proved 

\begin{fact}\label{Wilson}
For each countably infinite, $\omega$-categorical, characteristically simple group $H$, one of the following holds.
\begin{enumerate}
\item[(i)] $H$ is an elementary abelian $p$-group for some prime $p$.
\item[(ii)] $H \cong B(F)$ or $H \cong B^{-1}(F)$ for some non-abelian, finite, simple group $F$, where $B(F)$ is the group of all continuous functions from the Cantor set $\cal{C}$ to $F$, and $B^-(F)$ is the subgroup of $B(F)$ consisting of the functions $f$ such that $f(x_0)=e$ for a fixed element $x_0 \in \cal{C}$.
\item[(iii)] $H$ is a perfect $p$-group (perfect means that $H$ equals its commutator subgroup).
\end{enumerate}
\end{fact}
 
This suggests a method of proving structural results about certain $\omega$-categorical groups. Namely, sometimes induction on the maximal (finite) length of a series of characteristic subgroups allows one to reduce the situation to the case of characteristically simple groups and then to apply the above Wilson's theorem.

In \cite{5}, $\omega$-categorical groups and rings satisfying NIP are considered. It was proved in there that $\omega$-categorical rings with NIP are nilpotent-by-finite, and it was conjectured that $\omega$-categorical groups with NIP are nilpotent-by-finite, too.
This conjecture was shown to be true, but under the additional assumption that the group in question has fsg (finitely satisfiable generics). In fact, such groups are generically stable according to the terminology introduced in \cite[Section 6]{7}. It remains an open problem how to deal with $\omega$-categorical groups satisfying NIP in the situation when the fsg assumption is dropped.

In this paper, we do not drop the fsg assumption, but we do drop NIP. More precisely, we consider $\omega$-categorical groups which are generically stable (and some variants of this situation). Our main result says that each  $\omega$-categorical, generically stable group  is solvable-by-finite. The proof is by induction on the maximal length of a series of characteristic (in a generalized sense) subgroups. There are three main new ingredients in comparison to the proof of \cite[Theorem 3.4]{5}. First of all, because of the lack of NIP, we have to prove certain chain conditions for subgroups uniformly definable over Morley sequences in generically stable types. Next, using them, we eliminate items (ii) and (iii) form Wilson's theorem. This step is similar to the one in \cite{5}, but it requires more work,  e.g. an application of finite Ramsey theorem. Finally, in the last part of the proof, we have to deal with solvable iterated commutators instead of nilpotent commutators, which requires different computational ideas. 

The second author would like to thank Anand Pillay for suggesting that one of the lemmas in \cite{5} should be true without the NIP assumption, which was the starting point for considerations contained in this paper.

\section{Preliminaries}

Recall that a first order structure $M$ in a countable language is said to be $\omega$-categorical if, up to isomorphism, $Th(M)$ has at most one model of cardinality $\aleph_0$. By Ryll-Nardzewski's theorem, this is equivalent to the condition that for every natural number $n$ there are only finitely many $n$-types over $\emptyset$. Assume $M$ is $\omega$-categorical. If $M$ is countable or a monster model (i.e. a model which is $\kappa$-saturated and strongly $\kappa$-homogeneous for a big cardinal $\kappa$), two finite tuples have the same type over $\emptyset$ iff they lie in the same orbit under the action of the automorphism group of $M$, and hence for each natural number $n$, the automorphism group of $M$ has only finitely many orbits on $n$-tuples (which implies that $M$ is locally finite). Moreover, 
for any finite subset $A$ of such an $M$, a subset $D$ of $M$ is $A$-invariant iff $D$ is $A$-definable.

 Let $T$ be a first order theory. We work in a monster model $\C$ of $T$.

Let $p \in S(\C)$ be invariant over $A\subset \C$. We say that $(a_i)_{i \in \omega}$ is a Morley sequence in $p$ over $A$ if 
$a_i \models p |A a_{<i}$ for all $i$. Morley sequences  in $p$ over $A$ are indiscernible over $A$ and they have the same order type over $A$.
If $\C' \succ \C$ is a bigger monster model, then the generalized defining scheme of $p$ determines a unique $A$-invariant extension $\widetilde{p}\in S(\C')$ of $p$ (by the generalized defining scheme of $p$ we mean a family of sets $\{ p_i^\varphi: i \in I_\varphi\}$ (with $\varphi(x,y)$ ranging over all formulas without parameters) of complete types over $A$ such that $\varphi(x,c) \in p$ iff $c \in \bigcup_{i \in I_\varphi} p_i^{\varphi}(\C)$). By a Morley sequence in $p$ we mean a Morley sequence in $\widetilde{p}$ over $\C$. Finally, $p^{(k)}$ (where $k \in \omega \cup\{ \omega \}$) denotes the type over $\C$ of a Morley sequence in $p$ of length $k$.

\begin{definition}
A global type $p \in S(\C)$ is said to be generically stable if, for some small $A$, it is $A$-invariant and for each formula $\varphi(x;y)$, there is a natural number $m$ such that for any Morley sequence $(a_i: i< \omega)$ in $p$ over $A$ and any $b$ from $\C$, either less than $m$ $a_i$'s satisfy $\varphi(b;y)$ or less than $m$ $a_i$'s satisfy $\neg \varphi(b;y)$. 
In this definition, as a witness set $A$ one can take any (small) set over which $p$ is invariant. We will say that $p$ is generically stable over $A$ to express that $p$ is invariant over $A$ and generically stable.
\end{definition}

 Suppose $p \in S(\C)$ is $A$-invariant. Assuming NIP, there are various equivalent definitions of generic stability of $p$ (see \cite[Proposition 3.2]{7}). It turns out that in general (i.e. 
without the NIP assumption) all these definitions are implied by the definition given above (see \cite[Proposition 1]{8}). In particular, if $p$ is generically stable, then a Morley sequence in $p$ over $A$ is an indiscernible set over $A$. Some observations on these issues are contained Section 3.


%
%

\begin{proposition}\label{dcl}
Let $p=tp(a/\C)$ be a type generically stable over $A$, and assume that $b \in dcl(a)$. Then $tp(b/\C)$ is also generically stable over $A$.
\end{proposition}
{\em Proof.}
Let $\C' \succ \C$ be a bigger monster model containing $a$ and $b$. Let $g$ be a $\emptyset$-definable function such that $b=g(a)$.

First, we check that $tp(b/\C)$ is $A$-invariant. Consider any $f \in Aut(\C/A)$. We can extend it to an $\overline{f} \in Aut(\C'/A)$. Then, $tp(\overline{f}(a)/\C)=tp(a/\C)$ and $\overline{f}(b)=g(\overline{f}(a))$, so $tp(\overline{f}(b)/\C)=tp(b/\C)$. Thus, 
$tp(b/\C)$ is $A$-invariant.

Now, we check the main part of the definition of generic stability for $tp(b/A)$. Let $(h_i)_{i<\omega}$ be a Morley sequence in $p$ over $A$. We claim that $(g(h_i))_{i<\omega}$ is a Morley sequence in $tp(b/\C)$ over $A$. 
Indeed, for any formula $\varphi(x,g(h_0),\dots,g(h_n))$ in $tp(b/\C)|_{A,g(h_0),\dots,g(h_{n-1})}$, putting $\psi(x,y_0,\dots,y_{n-1})=\varphi(g(x),g(y_0),\dots,g(y_{n-1}))$ , 
we see that $\psi(x,h_0,\dots,h_{n-1}) \in p|_{A,h_0,\dots,h_{n-1}}$, so $\psi(h_n,h_0,\dots,h_{n-1})$ holds, and thus $\varphi(g(h_n),g(h_0),\dots,g(h_{n-1}))$ holds as well.

Consider any formula $\alpha(x;y)$. Choose $m$ as in the definition of generic stability for the type $p$ and the formula $\beta(x;y):=\alpha(x;g(y))$. Then, for any $b$, we have that either $| \{ i<\omega \models \alpha(b;g(h_i)\}| < m$ or $| \{ i<\omega \models \neg \alpha(b;g(h_i)\}| < m$.

The above observations show that $tp(b/\C)$ is generically stable over $A$. \hfill $\square$\\

Recall that a subset of a group $G$ is said to be left generic if finitely many left translates of this set cover $G$. A formula $\varphi(x)$ is left generic if the set $\varphi(G)$ is left generic. A type is said to be left generic if every formula in it is left generic.

\begin{definition}\label{fsg}
Let $G$ be a group definable in $\C$ by a formula $G(x)$.
$G$ has fsg (finitely satisfiable generics) if there is a global type $p$ containing $G(x)$ and a model $M\prec \C$, of cardinality less than the degree of saturation of $\mathfrak C$, such that for all $g$, $gp$ is finitely satisfiable in $M$.
\end{definition}

Let $G$ be a group definable in $\C$.
By $G^{00}$ we will denote the smallest type-definable subgroup of bounded index (if it exists). 
The following fact was proved in \cite[Section 4]{4}.

\begin{fact}\label{fsg}
Suppose $G$ has fsg, 
witnessed by $p$. Then:\\
(i) a formula is left generic iff it is right generic (so we say that it is generic),\\
(ii) $p$ is generic,\\
(iii) the family of nongeneric sets forms an ideal, so any partial generic type can be extended to a global generic type,\\
(iv) $G^{00}$ exists, it is type-definable over $\emptyset$, and it is the stabilizer of any global generic type of $G$.
\end{fact}

Recall \cite[Proposition 0.26]{3}.

\begin{fact}\label{uniqueness}
Suppose $G$ has fsg and $G^{00}$ is definable. Then $G^{00}$ has a unique global generic type.
\end{fact}

The next definition was introduced in \cite[Section 6]{7}.

\begin{definition}
$G$ is generically stable if it has fsg and some global generic type is generically stable.
\end{definition}

We say that a group definable in a non-saturated model is generically stable if the group defined by the same formula in a monster model is such. 

A few more observations and questions on generic stability are contained in Section 3.

We will say that $G$ is connected if it does not have a proper, definable subgroup of finite index; $G$ is absolutely connected if it does not have a proper, type-definable subgroup of bounded index.

\section{Main results}

The goal of this section is to prove the main result of the paper, namely Theorem \ref{main theorem}. As was mentioned in the introduction, first we will prove certain chain conditions. Then, we will use them to eliminate items (ii) and (iii) from Wilson's theorem for characteristically simple groups. Finally, using induction and certain computations involving commutators and centralizers, we will reduce the situation to characteristically simple groups, and we will be done. 

\begin{lemma}\label{chain condition}
Let $G$ be a group which is $\emptyset$-definable in $\mathfrak{C}$ by a formula $G(x)$. 
Assume that $p\in S_{1}(\mathfrak{C})\cap [G(x)]$
is generically stable over $A$.  \\
(i) Let $H(x,\bar{z};y)$ be a formula over $A$, defining a family of groups $H(G,c;g),g\in G, c\in G^{k}$, where $k=|\bar{z}|$.
Then there is some $n<\omega$ such that for any $c\in G^{k}$, 
$(g_{0},g_{1},\dots)\models p^{(\omega)}|_{Ac}$ and $i_{1}<\dots <i_{n}$, the following equality holds
\[\bigcap_{i<\omega}^{}H(G,c;g_{i})=H(G,c;g_{i_{1}})\cap \dots \cap H(G,c;g_{i_{n}}).\] \\
(ii) Let $H(x,y_{1},\dots,y_{k})$  be a formula over $A$, defining a family of groups \\ $H(G,h_{1},\dots ,h_{k}), h_{1},\dots ,h_{k} \in G$. Then for every $(g_{0},g_{1},\dots)\models p^{(\omega)}|_{A}$, there is some $n<\omega$ such that
$$\begin{array}{l}
\bigcap_{i_{1}<\dots <i_{k}}^{}H(G,g_{i_{1}},\dots ,g_{i_{k}})=\bigcap_{i_{1}<\dots <i_{k}<n}^{}H(G,g_{i_{1}},\dots, g_{i_{k}})= \\ 
=\bigcap_{i_{1}<\dots <i_{k}, \; i_{1},\dots ,i_{k} \in S}^{}H(G,g_{i_{1}},\dots ,g_{i_{k}}) 
\end{array}$$
for any set $S\subseteq \omega$ of cardinality $n$. 
\end{lemma}
{\em Proof.}
(i)
Let $m$ be such as in the definition of generic stability for $p$ and $H(x,\bar{z} ;y)$.
Fix any $c\in G^{k}$. Now, $H(x,y)$ will denote the formula $H(x,c;y)$.
Put $H_{i}=H(G,g_{i})$. We will show that $\bigcap_{i<\omega}^{}H_{i}=\bigcap_{i\leq 2m}^{}H_{i},$
which will complete the proof, due to the indiscernibility over $Ac$ of the set $\{g_{i}:i<\omega\}$.
Let $H=\bigcap_{i\leq m}^{}H_{i}$. 

Notice that $H=\bigcup_{m+1\leq i \leq 2m}^{}(H\cap H_{i})$ for if it is not the case, then any
$a\in H\backslash \bigcup_{m+1\leq i\leq 2m}^{}H_{i}$ contradicts the choice of $m$.
Now, there is some $j\in [m+1,2m]$ such that
\[(*) \bigcap_{m+1\leq i\leq 2m, i\neq j}(H\cap H_{i})\subseteq H_{j}.    \]
If not, pick some $a_{j}\in \bigcap_{m+1\leq i\leq 2m, i\neq j}(H\cap H_{i})\backslash H_{j}$ for every $j\in [m+1,2m]$. Then, $a_{m+1}a_{m+2}\dots a_{2m}\in H\backslash
\bigcup_{m+1\leq i\leq 2m}H_{i}$, a contradiction.

By $(*)$ and the indiscernibility over $Ac$ of the set $\{g_{i}:i<\omega \}$, we get that 
$$\bigcap_{m+1\leq i\leq 2m, i\neq j}(H\cap H_{i})\subseteq H_{l}$$ for all $l>2m$. 
Thus, $\bigcap_{i<\omega}^{}H_{i}=\bigcap_{i\leq 2m}^{}H_{i}$.\\
\newline
(ii) 
Put $H_{i_{1},\dots ,i_{k}}=H(G,g_{i_{1}},\dots , g_{i_{k}})$ for all $i_{1},\dots ,i_{k}$.  Fix any $i_{1}<\dots <i_{k}<\omega$. For any 
$1\leq j\leq k$, let $n_{j}$ be as in the conclusion of (i) for the following formula $\phi(x,y_{1},\dots , y_{i-1},y_{i+1},\dots , y_{k};y_{i}):=H(x,y_{1},\dots ,y_{k})$. Put $n=\max (n_{1},\dots ,n_{k})$ and $l=1+nk+\max(i_{1},\dots ,i_{k})$. Since  
$\{ g_{i}:i<\omega \}$ is indiscernible over $A$, we have that $g_{i_{1}},g_{l},g_{l+1},\dots$ is a Morley sequence in $p$ over $A,g_{i_{2}},\dots ,g_{i_{k}}$. By the choice of $n$, we conclude that
$H_{i_{1},\dots ,i_{k}}\supseteq \bigcap_{j_{1}\in [l,l+n-1]}H_{j_{1},i_{2},\dots ,i_{k}}$.
Repeating this argument $k$ times, we get
$$H_{i_{1},\dots ,i_{k}}\supseteq \bigcap_{j_{1}\in [l,l+n-1],j_{2}\in [l+n,l+2n-1],\dots, j_{k}\in [l+(k-1)n,l+kn-1]}H_{j_{1},\dots ,j_{k}}.$$
Similarly, for $j_{1}\in [l,l+n-1],j_{2}\in [l+n,l+2n-1],\dots, j_{k}\in [l+(k-1)n,l+kn-1]$, we have
$$H_{j_{1},\dots ,j_{k}}\supseteq \bigcap_{s_{1}\in [0,n-1],s_{2}\in [n,2n-1],\dots, s_{k}\in [(k-1)n,kn-1]}H_{s_{1},\dots ,s_{k}}.$$ Hence, $H_{i_{1},\dots ,i_{k}}\supseteq \bigcap_{j_{1}<\dots <j_{k}<kn}H_{j_{1},\dots ,j_{k}}$. 
We conclude that $$\bigcap_{i_{1}<\dots <i_{k}}^{}H(G,g_{i_{1}},\dots ,g_{i_{k}})=\bigcap_{i_{1}<\dots <i_{k}<kn}^{}H(G,g_{i_{1}},\dots ,g_{i_{k}}),$$ which completes the proof, because the set $\{g_{i}:i<\omega\}$ is indiscernible over $A$. \hfill $\square$

\begin{lemma}\label{main lemma}
Let $G$ be a $\emptyset$-definable group (by a formula $G(x)$) in $\C \models T$, where $T$ is an $\omega$-categorical theory. Assume that $G_{1}\unlhd G$ is infinite, $\emptyset$-definable, and characteristically simple in $(G,\C)$, i.e. it has no non-trivial, proper subgroup which is invariant under conjugations by elements of $G$ and invariant under $Aut(\C)$.
Suppose $p\in S_{1}(\C)\cap [G(x)]$ is a type  generically stable over $\emptyset$ and such that  $(\forall A\subseteq G)(\forall a\models p_{|A})(\forall g\in G)(a^{g} \models p_{|A^{g}})$ (e.g. this holds when $G=\C$ is a pure group). Let $(h_{i})_{i<\omega}$ be a Morley sequence in $p$. 
Assume that there is $tp(d/\C) \in S(\C) \cap [G_1(x)]$ such that $d \ne e$ and $d=f(h_{0},\dots ,h_{k-1})$ for some function $f$, which is $\emptyset$-definable in the language $\{\cdot\}$, and for some $k<\omega$. Then $G_{1}$ is abelian.
\end{lemma}
{\em Proof.}  Let $(g_{i})_{i<\omega}$ be a Morley sequence in $p$ over $\emptyset$. Put  $H_{i_{1},\dots,i_{k}}=C_{G_{1}}(f(g_{i_{1}},\dots ,g_{i_{k}}))$ for all $i_{1},\dots ,i_{n}<\omega$, and $H=\bigcap_{i_{1}<\dots <i_{k}}H_{i_{1},\dots,i_{k}}$. From Lemma \ref{chain condition}(ii), we have that there is some $n<\omega$ such that for every $S\subseteq \omega$ of cardinality $n$,  $$H=\bigcap_{i_{1}<\dots <i_{k},i_{1},\dots ,i_{k}\in S}H_{i_{1},\dots,i_{k}}.$$ 

We will show that $H$ is invariant under $Aut(\C)$. Take any $h\in Aut(\C)$.
Put $a_{i}=h(g_{i})$, and choose a Morley sequence $(b_{i})_{i<\omega}$ in $p$ over $\{a_{i},g_{i}:i<\omega\}$. Notice that the sequences $(g_i : i<\omega)^\frown (b_i : i<\omega)$ and $(a_i : i<\omega)^\frown (b_i: i<\omega)$ are Morley sequences in $p$ over $\emptyset$, and thus they are indiscernible as sets. Therefore, 
$$\begin{array}{l}
H=\bigcap_{i_{1}<\dots <i_{k}<n}C_{G_1}(f(b_{i_{1}},\dots ,b_{i_{k}}))=\bigcap_{i_{1}<\dots <i_{k}<\omega}C_{G_{1}}(f(a_{i_{1}},\dots ,a_{i_{k}}))= \\ =\bigcap_{i_{1}<\dots <i_{k}<\omega}h[C_{G_{1}}(f(g_{i_{1}},\dots ,g_{i_{k}}))]=h[H],
\end{array}$$ and so $H$ is invariant under $Aut(\C)$. 

Next, we will show that $H$ is normal in $G$.
Take any $g\in G$.
We have 
$$\begin{array}{l}
H^{g}=\bigcap_{i_{1}<\dots <i_{k}}C_{G_{1}}(f(g_{i_{1}},\dots ,g_{i_{k}}))^{g}=\bigcap_{i_{1}<\dots <i_{k}}C_{G_{1}}(f(g_{i_{1}},\dots ,g_{i_{k}})^{g})= \\ =\bigcap_{i_{1}<\dots <i_{k}}C_{G_{1}}(f(g_{i_{1}}^{g},\dots ,g_{i_{k}}^{g}))
\end{array}$$
(the last equality holds, because $f$ is $\emptyset$-definable in the language $\{\cdot\}$).
Using the assumptions about $p$, we see that $(g_{i}^{g})_{i<\omega}$ is a Morley sequence in $p$ over $\emptyset$.
As above, it easily follows that $H=\bigcap_{i_{1}<\dots <i_{k}}C_{G_{1}}(f(g_{i_{1}}^{g},\dots ,g_{i_{k}}^{g}))=H^{g}$,
so $H$ is normal in $G$.

Now, we will show that $H\neq \{e\}$. It follows from the assumptions on $G_{1}$ that $G_1$ is a characteristically simple group. Take a countable $(M,\cdot)\prec (G_{1},\cdot)$. Then $M$ is also a characteristically simple group, so, by Fact \ref{Wilson}, $M$ is either a $p$-group or it is isomorphic to a group of the form $B(F)$ or $B^{-}(F)$. 

If $M$ is a $p$-group, then $G_1$ is also a $p$-group, so $\langle\{f(g_{i_{1}},\dots ,g_{i_{k}}):i_{1}<\dots <i_{k}<n\}\rangle$ is a finite $p$-group, hence it has non-trivial center. As $Z(\langle\{f(g_{i_{1}},\dots ,g_{i_{k}}):i_{1}<\dots i_{k}<n\}\rangle)\subseteq H$, $H$ is also non-trivial.

Now, consider the case when $M$ is of the form $B(F)$ (when $M=B^{-}(F)$, the argument is similar). Take any $y_{0}$ in the Cantor set. It is easy to see that if finitely many elements of $B(F)$ have the same value at a point from the Cantor set, then the intersection of their centralizers is non-trivial. By finite Ramsey theorem, there is a number $R<\omega$ such that for every $f_{1},\dots,f_{R}\in B(F)$, there exist $1\leq i_{1}<\dots <i_{n} \leq R$ such that for all $1\leq j_{1}<\dots <j_{k} \leq n$ the value of the function $f(f_{i_{j_{1}}},\dots ,f_{i_{j_{k}}})$ at $y_{0}$ is the same.
We conclude that $M$, and hence $G_1$, satisfies the following sentence
$$\forall x_{1},\dots ,x_{R}\bigvee_{1\leq i_{1}<\dots < i_{n}\leq R} \bigcap_{1\leq j_{1}<\dots < j_{k}\leq n}
C(f(x_{i_{j_{1}}},\dots ,x_{i_{j_{k}}})) \neq \{e\}.$$ Thus, by the choice of $n$, $H$ is non-trivial. 

From these observations, and from the characteristic simplicity in $(G,\C)$ of $G_{1}$, we conclude that $H=G_{1}$. Thus, $Z(G_{1}) \neq \{e\}$. But $Z(G_{1})$ is normal in $G$ and invariant under $Aut(\C)$, so $Z(G_{1}) =G_{1}$. Hence, $G_{1}$ is abelian. \hfill $\square$ \\

Now, we have all the tools in order to prove the main results of the paper.

\begin{theorem}
We work in a monster model $\C$ of an $\omega$-categorical theory.
Let $G$ be a $\emptyset$-definable group having a global generic type $p$ which is generically stable over $\emptyset$ and such that for every $\emptyset$-definable, normal subgroup $L$ of $G$, we have 
$(\forall A\subseteq G/L)(\forall a\models p_{|A})(\forall g\in G)(a^{g} \models p_{|A^{g}})$. Then $G$ is solvable-by-finite.
\end{theorem}
{\em Proof.} 
We will show that $G$ has a $\emptyset$-definable, solvable subgroup of finite index. Of course, we can assume that $G$ is infinite. The proof will be by induction on the greatest natural number $n$ for which there is a series $\{e\}=G_{0}< G_{1}< \dots < G_{n}=G$ of $\emptyset$-definable (in $\C$) normal subgroups of $G$. Notice that then $G_{k}/G_{k-1}$ is characteristically simple in $(G, \C)$  for every $k\in\{1,\dots,n\}$.
 
If $n=1$, then by Lemma \ref{main lemma}, $G$ is abelian.
We turn to the induction step, where we assume that $n>1$. By induction hypothesis, there is a $\emptyset$-definable $H\unlhd G$ such that $[G:H]<\omega$ and $H/G_{1}$ is solvable (we leave to the reader checking that the group $G/G_1$ satisfies the hypothesis of the theorem; one should use here Proposition \ref{dcl}). So, in order to finish the proof, it is enough to show that $G_1$ is abelian.

Let $(g_{i})_{i< \omega}$ be a Morley sequence in $p$ over $\emptyset$. There exist $i<j<\omega$ such that  $g_{i}H=g_{j}H$. Then, $[g_{i},g_{j}]\in H$, so $[g_{i_{1}},g_{i_{2}}]\in H$ for all $i_{1},i_{2}<\omega$. Hence, from the solvability of $H/G_{1}$, we get that there is a minimal $k<\omega$ such that $\delta_{k}(g_{0},\dots ,g_{2^{k}-1})\in G_{1}$, where the iterated commutator $\delta_l$ is defined recursively as follows: 
$$\begin{array}{l}
\delta_0(a_1)=a_1, \\ 
\delta_{l+1}(a_1,\dots,a_{2^{l+1}})=[\delta_l(a_1,\dots,a_{2^{l}}),\delta_l(a_{2^l+1},\dots,a_{2^{l+1}})].
\end{array}$$

Notice first that we can assume that $k>0$. Indeed, if $k=0$, then $p \in S(\C) \cap [G_1(x)]$, so $G_1$ is abelian by Lemma \ref{main lemma}.\\[2mm]
{\bf Case 1.} $\delta_{k}(g_{0},\dots ,g_{2^{k}-1})=e$.\\[2mm]
Put $K=\bigcap_{i_{1}<\dots <i_{2^{k-1}}<\omega}C(\delta_{k-1}(g_{i_{1}},\dots ,g_{i_{2^{k-1}}}))$.
As in the proof of Lemma \ref{main lemma}, using Lemma \ref{chain condition}(ii), one can show that $K$ is a $\emptyset$-invariant, normal subgroup of $G$. Let us show now that $Z(K)$ is non-trivial. By Lemma \ref{chain condition}(ii), there is some $m<\omega$ for which $K=\bigcap_{i_{1}<\dots <i_{2^{k-1}}<m}C(\delta_{k-1}(g_{i_{1}},\dots ,g_{i_{2^{k-1}}}))$. Hence, by the assumption of Case 1, we get that $\delta_{k-1}(g_{m},\dots ,g_{m+2^{k-1}-1}) \in K$. On the other hand,  it is clear that $K \subseteq C(\delta_{k-1}(g_{m},\dots ,g_{m+2^{k-1}-1}))$. Thus,  $\delta_{k-1}(g_{m},\dots ,g_{m+2^{k-1}-1}) \in Z(K)$, and, by the choice of $k$, $\delta_{k-1}(g_{m},\dots ,g_{m+2^{k-1}-1}) \ne e$.

Summarizing, $Z(K)$ is a non-trivial, $\emptyset$-invariant, normal, abelian subgroup of $G$. So, by induction hypothesis, $G/Z(K)$ has a $\emptyset$-definable, solvable subgroup of finite index, hence so does $G$.\\[2mm] 
{\bf Case 2.} $\delta_{k}(g_{0},\dots ,g_{2^{k}-1}) \neq e$.\\[2mm]
If $G_{1}$ is infinite, then by Lemma \ref{main lemma}, we get that $G_{1}$ is abelian, so $H$ is solvable, and we are done. So, we may assume that $G_{1}$ is finite. Then, $[G:C(G_{1})]<\omega$, and hence $[g_{i_{1}},g_{i_{2}}]\in C(G_{1})$ for every $i_{1},i_{2}<\omega$.
Since $k>0$, we see that $\delta_{k}(g_{0},\dots ,g_{2^{k}-1})\in G_{1} \cap C(G_1)=Z(G_1)$, and so $Z(G_1)$ is non-trivial. But, $Z(G_1)$ is a $\emptyset$-invariant, normal subgroup of $G$ contained in $G_1$. Therefore, $G_1=Z(G_1)$, i.e. $G_1$ is abelian, and we are done.  \hfill $\square$\\


The following corollary immediately follows from the last theorem.

\begin{corollary}\label{cor}
(i) We work in a monster model $\C$ of an $\omega$-categorical theory.
Let $G$ be a $\emptyset$-definable group having a global generic type which is generically stable over $\emptyset$. Assume that each inner automorphism of $G$ is induced by an automorphism of $\C$. Then $G$ is solvable-by-finite.\\
(ii) Assume that $G$ is an $\omega$-categorical, pure group possessing a global generic type which is generically stable over $\emptyset$. Then $G$ is solvable-by-finite.
\end{corollary}

\begin{theorem}\label{main theorem}
Every $\omega$-categorical, generically stable group is solvable-by-finite.
\end{theorem}
{\em Proof.} 
By Corollary \ref{cor}(ii), it is enough to reduce the situation to the case of a pure group having a generic type which is generically stable over $\emptyset$.

Let $G$ be an infinite, generically stable group, definable in a monster model $\C$ of an $\omega$-categorical theory. 
Let $p$ be a generically stable (global) type in $G$. By Fact \ref{fsg}(iv), $G^{00}$ exists. Since $G^{00}$ is $\emptyset$-invariant, $\omega$-categoricity implies that it is $\emptyset$-definable. Hence, $[G:G^{00}]<\omega$. By Fact \ref{uniqueness}, $G^{00}$ has a unique generic type $p'$, which is a translation of $p$, so it is also generically stable (by Proposition \ref{dcl}) and witnesses fsg in $G^{00}$. Let $M\prec \C$ be a small model in which $p'$ is finitely satisfiable. Let $q\in S(G^{00})$ be a type naturally determined by $p'$ (we consider $G^{00}$ with the structure induced from $\C$). It is easy to check that $M\cap G^{00}\prec G^{00}$, and that $q$ is finitely satisfiable in $M\cap G^{00}$ and invariant under translations by elements of $G^{00}$ (by Fact \ref{fsg}(iv)). Put $r=q |_{\{\cdot\}}$. Then $r$ witnesses fsg in $(G^{00},\cdot)$, and $(G^{00},\cdot)$ is absolutely connected  (notice that we can choose $\C$ so that $(G^{00},\cdot)$ is a monster model of its theory). Hence, by Fact \ref{uniqueness}, $r$ is a unique generic type in  $(G^{00},\cdot)$. Thus, it is $\emptyset$-invariant and also generically stable over $\emptyset$, because a Morley sequence in $p'$ over $\emptyset$ is also a Morley sequence in $q$ over $\emptyset$, which is a Morley sequence in $r$ over $\emptyset$ (and there is only one Morley sequence up to the type). 
Hence, replacing $G$ by $(G^{00},\cdot)$, we can assume that $G$ is a pure group, $G=G^{00}$ and $G$ has a generic type which is  generically stable over $\emptyset$. 
 \hfill $\square$

\section{Remarks and questions about generic stability}

We work in a monster model $\C$ of a theory $T$.
In \cite{7}, among others, the following properties of a global, $A$-invariant type $p$ are considered:

\begin{enumerate}
\item[(i)] $p$ is generically stable over $A$,
\item[(ii)] $p$ is definable and finitely satisfiable in some small model containing $A$,
\item[(iii)] a Morley sequence in $p$ over $A$ is indiscernible over $A$ as a set.
\end{enumerate}

Under the NIP assumption, it is proved in \cite{7} that all these properties are equivalent.
Moreover, the proofs of implications $(i)\Rightarrow (ii)$ and $(ii)\Rightarrow (iii)$ do not use the NIP assumption, so these implications are true also without  NIP. Below we give an example showing that without NIP the implication $(iii)\Rightarrow (ii)$ is not always true. As to the implication $(ii) \Rightarrow (i)$, we think it is also false in general, but we have not found an appropriate example.\\[5mm]
%
%
{\bf Example} Let $(\C,R)$ be a monster model of the theory of random graphs. This theory is complete, has quantifier elimination, and is even $\omega$-categorical. Let $p$ be the global type determined by the collection of formulas $\{R(x,a)\land x\neq a:a\in \C\}$. We see that $p$ is $\emptyset$-invariant. Moreover, if $(g_i: i \in \omega)$ is a Morley sequence in $p$ over $\emptyset$, then, since $R$ is symmetric, $R(g_i,g_j)$ holds for all pairwise distinct $i,j<\omega$, and so $(g_i : i \in \omega)$ is indiscernible as a set (because we have quantifier elimination). On the other hand, if $M \prec \C$ is a small model, then the set of formulas $\{ \neg R(x,b) : b \in M \}$ is consistent, and for $a \in \C$ realizing all formulas from this set,  the formula $R(x,a)$ belongs to $p$ but has no realization in $M$. This shows that (iii) does not imply (ii).


\begin{remark}
Let $G=\C$ be an $\omega$-categorical group with fsg. Then there is a generic type in $G$ which is 
definable over $\emptyset$ and finitely satisfiable in a small model $M\prec G$. 
\end{remark}
{\em Proof.} By $\omega$-categoricity, $G^{00}$ is $\emptyset$-definable. 
From Fact \ref{uniqueness}, we get that in $G^{00}$ there is a unique global generic type $p$. 
So, $p$ is $\emptyset$-invariant, and, by $\omega$-categoricity,  it is $\emptyset$-definable. Moreover, by fsg, $p$ is finitely satisfiable in a small model $M\prec G$. \hfill $\square$ \\

By this remark, we see that if $G$ is an $\omega$-categorical group with fsg, then it has a 
$\emptyset$-invariant generic type satisfying property (ii). 

\begin{question}
Is it true, that if $G$ is an $\omega$-categorical group with fsg, then $G$ has a generically stable (over 
$\emptyset$) generic type?
\end{question}

If the answer is affirmative, then in Theorem \ref{main theorem} it is enough to assume that $G$ is $\omega$-categorical and has fsg.

We also ask the following question about generically stable types in arbitrary theories.

\begin{question}
Is it true, that if a global type $p$ is generically stable, then for every $n<\omega$, the type $p^{(n)}$ is also generically stable?
(Recall that $p^{(n)}=tp(g_{0},\dots ,g_{n-1}/\C)$, where $(g_{i})_{i<n}$ is a Morley sequence in $p$.)
\end{question}

An affirmative answer would allow us to simplify the proofs of some results of Section 2. Namely, Lemma \ref{chain condition}(ii) could be easily deduced from Lemma \ref{chain condition}(i), and in Lemma \ref{main lemma} we would get from the assumptions that $tp(d/\C)$ is generically stable, which would slightly simplify the proof of this lemma.

We finish with a few questions about $\omega$-categorical rings. In \cite{5}, it was proved that $\omega$-categorical rings with NIP are nilpotent-by-finite. One can ask what can be said about $\omega$-categorical rings when the NIP assumption is replaced by the fsg assumption for the additive group and/or by the assumption of  the existence of an additive generic type which is generically stable over $\emptyset$. For example, we have

\begin{question}\label{last question}
Are $\omega$-categorical, generically stable rings nilpotent-by-finite?
\end{question}

The fact that $\omega$-categorical rings with NIP are nilpotent-by-finite was used in \cite{5} to show that $\omega$-categorical groups with NIP and fsg are nilpotent-by-finite. Do we have something like that without the NIP assumption? More precisely, 
if the answer to Question \ref{last question} is positive, does it help to strengthen the conclusion of Theorem \ref{main theorem} by saying that the group in question is nilpotent-by-finite?

\noindent
{\bf Address:}\\
Instytut Matematyczny, Uniwersytet Wroc\l awski,\\
pl. Grunwaldzki 2/4, 50-384 Wroc\l aw, Poland.\\[3mm]
{\bf E-mail addresses:}\\
Jan Dobrowolski: Jan.Dobrowolski@math.uni.wroc.pl \\
Krzysztof Krupi\'nski:  kkrup@math.uni.wroc.pl


\begin{thebibliography}{99}

\bibitem{1} W. Baur, G. Cherlin, A. Macintyre. 
{\em Totally categorical groups and rings}, Journal of Algebra
57, 407-440, 1979. 


\bibitem{3} C. Ealy, K. Krupi\'nski, A. Pillay,
{\em Superrosy dependent groups having finitely satisfiable generics}, Annals of Pure and Applied Logic 151,
1-21, 2008.

\bibitem{EW}
D. Evans, F. Wagner, {\em Supersimple $\omega$-categorical groups and theories}, Journal of Symbolic Logic 65, 767-776, 2000.

\bibitem{7} E. Hrushovski, A. Pillay, 
{\em On NIP and invariant measures}, Journal of the European Mathematical Society, accepted.

\bibitem{4} E. Hrushovski, A. Pillay, Y. Peterzil,
{\em Groups, measures, and the NIP}, Journal of the American Mathematical Society 21, 563-595, 2008. 

\bibitem{5} K. Krupi\'nski,
{\em On $\omega$-categorical groups and rings with NIP}, Proceedings of the American Mathematical Society, accepted.

\bibitem{6} H. D. Macpherson, {\em Absolutely ubiquitous structures and $\aleph_{0}$-categorical groups}, Quart. J.
Math. Oxford (2) 39, 483-500, 1988. 

\bibitem{8} A. Pillay, P. Tanovi\'c, {\em Generic stability, regularity, and quasi-minimality}, preprint, 2009.

\bibitem{9} B. Poizat,
{\em  Stable groups}, American Mathematical Society, Providence, 2001.

\bibitem{10} J. Wilson,
{\em The algebraic structure of $\omega$-categorical groups}, in: Groups-St. Andrews, Ed. C.
M. Campbell, E. F. Robertson, London Math. Soc. Lecture Notes 71, Cambridge, 345-358,
1981. 



\end{thebibliography}
\end{document}